\documentclass[11pt]{amsart}
\usepackage{amssymb}
\usepackage{url}
\usepackage{enumerate}
\usepackage{tikz}
\newtheorem{lemma}{\bf Lemma}
\newtheorem{theorem}[lemma]{\bf Theorem}
\newtheorem{corollary}[lemma]{\bf Corollary}

\newtheorem{example}[lemma]{\bf Warm-Up Example}

\newtheorem{conjecture}[lemma]{\bf Conjecture}
\newtheorem{notation}[lemma]{\bf Notation}

\DeclareMathOperator{\id}{id}
\DeclareMathOperator{\dl}{dl}
\DeclareMathOperator{\st}{st}
\DeclareMathOperator{\lk}{lk}
\DeclareMathOperator{\BW}{BW}
\DeclareMathOperator{\BWI}{BWI}

\begin{document}
\parskip = 0mm
\title[A Conjecture of Kozlov from 1998]{A Conjecture of Kozlov from the 1998 {\sl Proceedings of the American Mathematical Society}\\ \footnotesize Non-Evasive Order Complexes and Generalizations of Non-Complemented Lattices}
\author{Jonathan David Farley}
\address{Department of Mathematics, Morgan State University, 1700 E. Cold Spring Lane, Baltimore, MD 21251, United States of America, {\tt lattice.theory@gmail.com}}

\keywords{(Partially) ordered set, simplicial complex, order complex, non-evasive, (complemented) lattice, chain.}

\subjclass[2010]{05E45, 06A07, 06B75, 06C15}

\begin{abstract}
Let $P$ be a finite poset with an element $s$ such that
\begin{itemize} \item[(1)] for all $x\in P$, either $s\vee x$ or $s\wedge x$ exists; and
\item[(2)] for all $x,y\in P$ such that $x<y$, if $s\wedge x$ does not exist but $s\wedge y$ does exist, then $(s\wedge y)\vee x$ exists.
\end{itemize}

Kozlov conjectured in the 1998 {\it Proceedings of the American Mathematical Society} that the order complex of $P$ is non-evasive.

We prove this conjecture.
\end{abstract}

\thanks {The author would like to thank Dr. Marshall Cohen for comments regarding this paper.}

\maketitle


\def\Qa{\mathbb{Q}_0}
\def\Qb{\mathbb{Q}_1}
\def\Q{\mathbb{Q}}
\def\card{{\rm card}}
\parskip = 2mm
\parindent = 10mm
\def\Part{{\rm Part}}
\def\P{{\mathcal P}}
\def\Eq{{\rm Eq}}
\def\cld{Cl_\tau(\Delta)}
\def\Csing{{\mathcal C}_{\{*\}}}
\def\Cftwo{{\mathcal C}_{{\rm fin}\rangle1}}
\def\Cinf{{\mathcal C}_{\infty}}
\def\Pcf{{\mathcal P}_{\rm cf}}
\def\Fn{{\mathcal F}_n}
\def\proof{{\it Proof. }}


\vspace*{-4mm} 


In the 1998 {\sl Proceedings of the American Mathematical Society}, Kozlov, the winner of the 2005 European Prize in Combinatorics (``for deep combinatorial results obtained by algebraic topology and particularly for the solution of a conjecture of Lov\'asz'' \cite{FelLueNesJG}), made the following conjecture \cite[Conjecture 2.6]{KozIH}:

\begin{conjecture}[] Let $P$ be a finite poset with an element $s$ such that
\begin{itemize} \item[(1)] for all $x\in P$, either $s\vee x$ or $s\wedge x$ exists; and
\item[(2)] for all $x,y\in P$ such that $x<y$, if $s\wedge x$ does not exist but $s\wedge y$ does exist, then $(s\wedge y)\vee x$ exists.
\end{itemize}

Then the order complex of $P$ is non-evasive.
\end{conjecture}
(Definitions follow.)

We prove this conjecture (Corollary 15).

For the order theory background, terminology and notation, see \cite{DavPriJB} and \cite{SchAF} and their references.  A subscript (e.g., $a\wedge_Q b$) indicates the poset or subposet with respect to which the order-theoretic construction, relation, or operation is being considered.

Let $P$ be a poset.  If $x,y\in P$, we say $x$ is a {\it lower cover} of $y$ and $y$ is an {\it upper cover} of $x$ if $x<y$ and $\uparrow x\cap\downarrow y=\{x,y\}$, where for $z\in P$, $\downarrow z=\{w\in P\mid w\le z\}$ and $\uparrow z=\{w\in P\mid z\le w\}$.  For $Q\subseteq P$, $Q^{\it u}=\cap_{q\in Q} \uparrow q$ and $Q^\ell=\cap_{q\in Q} \downarrow q$.  We say $x\in P$ is {\it irreducible} if it has a unique lower cover or a unique upper cover.  The non-empty $n$-element poset $P$ is {\it dismantlable by irreducibles} if $P=\{p_1,\dots,p_n\}$ and $p_i$ is irreducible in the subposet $\{p_1,\dots,p_i\}$ for $i=2,3,\dots,n$.

If $f,g:P\to P$ are order-preserving maps, $f\le g$ means that $f(p)\le g(p)$ for all $p\in P$; denote the poset of order-preserving maps by $P^P$ and denote the constant map to $p\in P$ by $\langle p\rangle$.  Theorem 4.1 of \cite{BacBjoGI} says:

\begin{theorem} A finite non-empty poset $P$ is dismantlable by irreducibles if and only if the identity map $\id_P$ and some constant map are in the same connected component of $P^P$. $\qed$
\end{theorem}

A {\it simplicial complex} is a down-set $\Sigma$ of a power set of a set.  (Deviating from \cite[\S9]{BjoIE}, we will allow $\emptyset$ to be a member of $\Sigma$; but note that \cite{KozIH} refers to \cite[9.9]{BjoIE}.)  The set of {\it vertices} of $\Sigma$ is $V:=V(\Sigma):=\big\{x\mid\{x\}\in\Sigma\big\}$.  If $W\subseteq V$, $\Sigma|W$ denotes the simplicial complex $\{\sigma\in\Sigma\mid\sigma\subseteq W\}$.  For $\sigma\in\Sigma$, $\dl_\Sigma(\sigma)$ is $\Sigma|(V\setminus\sigma)$; $\st_\Sigma(\sigma)$ is $\{\tau\in\Sigma\mid\sigma\cup\tau\in\Sigma\}$; $\lk_\Sigma(\sigma)$ is $\dl_\Sigma(\sigma)\cap\st_\Sigma(\sigma)$; if $v\in V$, we write $\st_\Sigma(v)$ instead of $\st_\Sigma(\{v\})$, and $\Sigma$ is a {\it cone} with {\it peak} $v$ if $\Sigma=\st_\Sigma(v)$.

We define the term ``non-evasive'' by induction on $|V|$.  For $V$ finite, $\Sigma$ is {\it non-evasive} if, for some $v\in V$, $\Sigma=\big\{\emptyset,\{v\}\big\}$ or if $|V|>1$ and both $\dl_\Sigma(v)$ and $\lk_\Sigma(v)$ are non-evasive.  The term stems from an interesting connection with algorithms.  See the proof of \cite[Proposition 1]{KahSakStuHD} for that and for the connection with the notion of collapsibility, and see \cite[Definition 13]{BacAB}, which uses the term {\it link reducible}.

A simplicial complex whose realization is a set of two points on the real line---abstractly, $\Sigma=\big\{\emptyset,\{0\},\{1\}\big\}$---is not non-evasive, as $\lk_\Sigma(0)=\{\emptyset\}$, which is not non-evasive, having no vertices. But a simplicial complex whose realization is the unit interval---abstractly, $\Delta=\big\{\emptyset,\{0\},\{1\},\{0,1\}\big\}$---is non-evasive, as $|V|=2$ and $\dl_\Delta(0)=\lk_\Delta(0)=\big\{\emptyset,\{1\}\big\}$, which is non-evasive.

We will use \cite[Corollary 14]{BacAB}:

\begin{theorem} Let $\Sigma$ be a non-empty simplicial complex with $V=V(\Sigma)$ finite. Let $W\subseteq V$ and let $U=V\setminus W$.  If $\st_\Sigma(\sigma)|U$ is non-evasive for all $\sigma\in\Sigma|W$, then $\Sigma$ is non-evasive. $\qed$
\end{theorem}

If $P$ is a poset, the order complex $\Delta(P)$ is the family of chains of $P$.

The conjecture of Kozlov was inspired by \cite[Theorem 3.2]{BjoWalHC}, stating that such posets are contractible, and by Kozlov's own proof that the order complexes of non-empty finite truncated lattices with certain elements removed are non-evasive \cite[Theorem 2.4]{KozIH}.

The following is stated on \cite[p. 3462]{KozIH} and in \cite[Corollary 25]{BacAB}:

\begin{example} A finite cone is non-evasive.
\end{example}

\proof Let $\Sigma$ be a cone---say, $\Sigma=\st_\Sigma(v)$ where $v\in V=V(\Sigma)$.  If $\Sigma=\big\{\emptyset,\{v\}\big\}$, then $\Sigma$ is non-evasive.

So assume there exists $x\in V\setminus\{v\}$.

\quad {\it Claim 1. The complex $\dl_\Sigma(x)$ is a cone.}

{\it Proof of claim.} Note that $v\in V\big(\dl_\Sigma(x)\big)$ and $\st_{\dl_\Sigma(x)}(v)=\dl_\Sigma(x)$. $\qed$

\quad {\it Claim 2. The complex $\lk_\Sigma(x)$ is a cone.}

{\it Proof of claim.} We show that $\lk_\Sigma(x)=\st_{\lk_\Sigma(x)}(v)$.  Note that $v\in V\big(\lk_\Sigma(x)\big)$.

Let $\sigma\in\lk_\Sigma(x)$.  Then $\sigma\cup\{x\}\in\Sigma$ but $x\notin\sigma$.  Thus, since $\Sigma=\st_\Sigma(v)$, $\sigma\cup\{v,x\}\in\Sigma$, so $\sigma\cup\{v\}\in\Sigma$ and $x\notin\sigma\cup\{v\}$.  Hence $\sigma\cup\{v\}\in\lk_\Sigma(x)$ and $\sigma\in\st_{\lk_\Sigma(x)}(v)$. $\qed$

Since $x\notin V\big(\dl_\Sigma(x)\big),V\big(\lk_\Sigma(x)\big)$, then by induction on $|V|$ Claims 1 and 2 imply that $\dl_\Sigma(x)$ and $\lk_\Sigma(x)$ are non-evasive.  Hence $\Sigma$ is non-evasive. $\qed$

Example 6 comes from \cite[Proposition 23 and Corollary 25]{BacAB}, but we can prove it more directly, so why not?

\begin{lemma} Let $P$ be a finite poset with an element $x$ that has a unique lower cover $x_*$.  Then $\lk_{\Delta(P)}(x)$ is a cone with peak $x_*$.  Hence, if $\dl_{\Delta(P)}(x)$ is non-evasive, then $\Delta(P)$ is non-evasive.
\end{lemma}

\proof Now $\lk_{\Delta(P)}(x)$ consists of all the chains $C=\{c_1,\dots,c_n\}$ of $P$ such that $x\notin C$ but $C\cup\{x\}$ is a chain.  Assume $c_1<c_2<\cdots<c_{k-1}<x<c_k<\cdots<c_n$ for some $n\in\mathbb N_0$ and $k\in\{1,\dots, n,n+1\}$.  That is, $x_*\in V\big(\lk_{\Delta(P)}(x)\big)$ and $\lk_{\Delta(P)}(x)=\st_{\lk_{\Delta(P)}(x)}(x_*)$ because $c_{k-1}\le x_*<x$, so $\lk_{\Delta(P)}(x)$ is a cone with peak $x_*$.

Now use Example 4. $\qed$

\begin{example} If $P$ is a finite dismantlable poset, then $\Delta(P)$ is non-evasive.
\end{example}

\proof The proof is by induction on $|P|$.

 \quad {\it Base Case. $|P|=1$}

Then $\Delta(P)$ is non-evasive!

\quad {\it Induction Step.}

Without loss of generality, assume $x\in P$ has a unique lower cover $x_*$.  Then $\Delta(P\setminus\{x\})=\dl_{\Delta(P)}(x)$ is dismantlable, hence non-evasive.

So, by Lemma 5, $\Delta(P)$ is non-evasive. $\qed$

The statement of Theorem 8 is modeled after \cite[Conjecture 2.6]{KozIH}---but {\sl vorsicht!}---and its proof mimics that of \cite[Theorem 26]{BacAB}.

\begin{lemma} Let $P$ be a poset.  Let $a,b,c\in P$ be such that $a\wedge b$ and $(a\wedge b)\wedge c$ exist.  Then $\bigwedge\{a,b,c\}$ exists and equals $(a\wedge b)\wedge c$.

If also $b\wedge c$ exists, then $\bigwedge\{a,b,c\}=a\wedge(b\wedge c)$.
\end{lemma}

\proof Now $a,b\ge a\wedge b$ so $a,b,c\ge(a\wedge b)\wedge c$.  If $x\in\{a,b,c\}^\ell$, then $x\in\{a,b\}^\ell$, so $a\wedge b\ge x$, implying that $x\in\{a\wedge b,c\}^\ell$ and thus $(a\wedge b)\wedge c\ge x$.  Hence $(a\wedge b)\wedge c=\bigwedge\{a,b,c\}$.

Now assume also that $b\wedge c$ exists.  We have $\bigwedge\{a,b,c\}\le a,b\wedge c$.  If $x\in\{a,b\wedge c\}^\ell$, then $x\le a,b,c$, so $x\le\bigwedge\{a,b,c\}$.

Hence $\bigwedge\{a,b,c\}=a\wedge(b\wedge c)$. $\qed$

\begin{theorem} Let $P$ be a finite poset.  Let $s\in P$.  

Assume that for each $y\in P$, either $s\wedge y$ exists or $s\vee y$ exists.

Assume that if $w,z\in P$ and $w<z$ and $s\wedge w$ does not exist but $s\wedge z$ does exist, then $z\wedge(w\vee s)$ exists.

Then $\Delta(P)$ is non-evasive.
\end{theorem}

\proof Let $\Sigma=\Delta(P)$.  Let
$$
W=\{y\in P\mid s\wedge y \text{ does not exist}\}.
$$
\noindent Let $U=P\setminus W$.  Thus $P=U\cup W$ and $U\cap W=\emptyset$.  To apply Theorem 3, we must show that $\st_\Sigma(\sigma)|U$ is non-evasive for all $\sigma\in\Sigma|W=\Delta(W)$.

First let $\sigma=\emptyset$.  Now $[\st_\Sigma(\emptyset)]|U=\Delta(U)$.  Define $f:U\to U$ as follows: for $z\in U$, let $f(z)=s\wedge z$.  Note that $s\ge s\wedge z$, so $s\wedge(s\wedge z)$ exists; hence $f(z)=s\wedge z\in U$.  Clearly $f\in U^U$.  Also, $s\in U$.

Note that $\id_U\ge f\le\langle s\rangle$ in $U^U$, so, by Theorem 2, $U$ is dismantlable by irreducibles.  Example 6 says $\Delta(U)$ is non-evasive.

Now assume $\sigma$ is a non-empty simplex of $\Delta(W)$---say, $\sigma=\{w_1,\dots,w_n\}$ where $n\in\mathbb N$ and $w_1<w_2<\cdots<w_n$.  Then $\st_\Sigma(\sigma)|U$ consists of those chains $\tau$ of $U$ such that $\sigma\cup\tau$ is a chain of $P$---say $$
\tau=\{u_{01},\dots,u_{0m_0},u_{11},\dots,u_{1m_1},\dots,u_{n1},\dots,u_{nm_n}\}
$$
\noindent where
$$
u_{01}<\cdots<u_{0m_0}<w_1<u_{11}<\cdots<u_{1m_1}<w_2<\cdots<w_n<u_{n1}<\cdots<u_{nm_n}
$$
\noindent with $m_0,m_1,\dots,m_n\in\mathbb N_0$.

Let
$$
T=U\cap[\downarrow w_1\cup(\uparrow w_1\cap\downarrow w_2)\cup(\uparrow w_2\cap\downarrow w_3)\cup\cdots\cup(\uparrow w_{n-1}\cap\downarrow w_n)\cup\uparrow w_n].
$$
\noindent  So $\st_\Sigma(\sigma)|U=\Delta(T)$.

To finish the proof of the theorem, it suffices to show that $T$ is dismantlable by irreducibles and hence $\Delta(T)$ is non-evasive by Example 6.

Note that any $t\in T$ is in a unique one of the sets in the above union, since $U\cap W=\emptyset$.

Let $i\in\{1,\dots,n\}$.  Note that $s\wedge w_i$ does not exist, so $s\vee w_i$ exists.  Also, $s\le s\vee w_i$, so $s\vee w_i\in U$.  In particular, $s\vee w_n\in U\cap \uparrow w_n$, so $T\ne\emptyset$.

Assume $t\le w_1$.  Then $f(t)=s\wedge t$ exists and is in $U$.  Since it is in $U\cap\downarrow w_1$, it is in $T$.

Assume $w_i\le t$.  Again, $s\wedge t$ exists. As $w_i<t$, our hypothesis tells us $t\wedge(w_i\vee s)$ exists.

Obviously $s\wedge(w_i\vee s)=s$ exists and we already know that $f(t)=t\wedge s=t\wedge[(w_i\vee s)\wedge s]$ exists so, by Lemma 7, $[t\wedge(w_i\vee s)]\wedge s$ exists (and equals $f(t)$, which is in $U$).  Thus $t\wedge(w_i\vee s)\in U$.  Also $w_i\le t\wedge(w_i\vee s)\le t$.  If also $t\le w_{i'}$ for $i'\in\{1,\dots,n\}$, then $t\wedge(w_i\vee s)\le w_{i'}$.

Now define $g:T\to T$ as follows: For $t\in T$,
$$
g(t)=
\begin{cases}
t\wedge s \ \text{if}\ t\le w_1, \\
t\wedge(w_i\vee s)\ \text{if}\ i\in\{1,\dots,n\}\text{ is the greatest number}\ k\ \text{ such that }\ w_k\le t.
\end{cases}
$$
\noindent We see that $g\in T^T$ and
$$
\id_T\ge g\le\langle w_n\vee s\rangle\text{,}
$$
\noindent proving that $T$ is dismantlable by Theorem 2. $\qed$

Wait---{\sl that's} not the conjecture we originally wanted to prove!  Let's try again.

\begin{notation} Let $\BW(P,s)$ be the statement:

``Let $P$ be a finite poset with an element $s\in P$.  Assume that for all $a\in P$, $a\vee s$ exists or $a\wedge s$ exists.  Assume that for all $a,b\in P$ if $a>b$ and $a\vee s$ does not exist and $b\vee s$ does exist, then $a\wedge(b\vee s)$ exists.''

Let $\BW(P,s,r)$ be the statement:

``$\BW(P,s)$ holds and $r$ is a minimal element of $P\setminus(\uparrow s\cup\downarrow s)$.''

Let $\BWI(P,s,r)$ be the statement:

``$\BW(P,s,r)$ holds and $r$ does not have a unique lower cover.''
\end{notation}

\begin{corollary} Let $P$ be a finite poset of the form $\uparrow s\cup\downarrow s$ for some $s\in P$.  Then $\Delta(P)$ is a cone with peak $s$.  Hence, $\Delta(P)$ is non-evasive.
\end{corollary}

\proof Use Example 4. $\qed$

\begin{lemma} Let $\BW(P,s,r)$ hold.  Then each lower cover of $r$ is in
$$
(\downarrow s)\setminus\{s\}. \qed
$$
\end{lemma}

\begin{lemma} Let $\BW(P,s,r)$ hold.  Then $\BW(P\setminus\{r\},s)$ holds.
\end{lemma}

\proof Since $r\notin\uparrow s\cup\downarrow s$, then $r\ne s$ so $s\in P\setminus\{r\}$.  Let $a\in P\setminus\{r\}$.  Assume that $a\vee s$ exists. Then $a\vee s\ne r$ (since $s\le a\vee s$ but $s \nleqslant r$) so $a\vee s\in P\setminus\{r\}$ and $a\vee s=a\vee_{P\setminus\{r\}}s$.

Now assume $a\wedge s$ exists.  Then $a\wedge s\ne r$ (since $s\ge a\wedge s$ but $s\ngeqslant r$) so $a\wedge s\in P\setminus\{r\}$ and $a\wedge s=a\wedge_{P\setminus\{r\}} s$.

Now assume that $a,b\in P\setminus\{r\}$, $a>b$, $a\vee_{P\setminus\{r\}} s$ does not exist and $b\vee_{P\setminus\{r\}} s$ does exist.  If $z\in\{b,s\}^{\it u}$ then $z\ne r$ (since $r\ngeqslant s$) so $b\vee_{P\setminus\{r\}} s\le z$ and hence $b\vee_{P\setminus\{r\}} s=b\vee s$.

We have already shown that $a\vee s$ does not exist (or else $a\vee_{P\setminus\{r\}} s$ would exist).  By $\BW(P,s)$ $a\wedge(b\vee s)$ exists.  But this is $a\wedge(b\vee_{P\setminus\{r\}} s)$.

{\it Case 1. $r=a\wedge(b\vee_{P\setminus\{r\}} s)$}

Then $r\le a$.  We know that $a\wedge s$ exists and is in $P\setminus\{r\}$ so $a\wedge s<r$.  Thus $r$ has a lower cover, $r_1$.  If $r_2$ is also a lower cover, then by Lemma 11 $r_1,r_2<s$ and $r_1,r_2<r\le a$, so $r_1,r_2\le a\wedge s<r$.  Hence $r_1=r_2=a\wedge s$ and $r$ has a unique lower cover.  If $z\in (P\setminus\{r\})\cap\{a, b\vee_{P\setminus\{r\}} s\}^\ell$, then $z<r$ so $z\le r_1$.  Hence $r_1= a\wedge_{P\setminus\{r\}} (b\vee_{P\setminus\{r\}} s)$.

{\it Case 2. $r\ne a\wedge(b\vee_{P\setminus\{r\}} s)$}

Then $a\wedge (b\vee_{P\setminus\{r\}} s)\in P\setminus\{r\}$ so it equals  $a\wedge_{P\setminus\{r\}} (b\vee_{P\setminus\{r\}} s)$.

This establishes $\BW(P\setminus\{r\},s)$. $\qed$

\begin{lemma} Let $\BWI(P,s,r)$ hold.  Then $\BW\big((\uparrow r\cup\downarrow r)\setminus\{r\},r\vee s\big)$ holds.
\end{lemma}

\proof Assume for a contradiction that $r\vee s$ does not exist.  Then by $\BW(P,s)$, $r\wedge s$ exists.  Since $r\nleqslant s$, we have $r\wedge s<r$.  Thus $r$ has a lower cover---by $\BWI(P,s,r)$, at least two, $r_1$ and $r_2$.  By Lemma 11, $r_1,r_2<s$ so $r_1,r_2\le r\wedge s<r$ and hence $r_1=r\wedge s=r_2$, a contradiction.

Obviously $r\le r\vee s$ and $r\ne r\vee s$ since $s\nleqslant r$.  Let $q:=r\vee s$.  Let $Q:=(\uparrow r\cup\downarrow r)\setminus\{r\}$.  

\quad {\it Claim. If $b\in Q$ and $b\vee s$ exists, then $b\vee_Q q$ exists.}

{\it Proof of claim.} If $b<r$, then as $q\ge r$, this means $b\vee q=q$ so $b\vee_Q q=q$.

Otherwise, $b>r$, so $b=b\vee r$.  We have $b\vee s=(b\vee r)\vee s$, which by Lemma 7 is $b\vee(r\vee s)=b\vee q$.  As $b\vee q\ge q>r$, we have $b\vee q\in Q$, so $b\vee q=b\vee_Q q$. $\qed$

Let $a\in Q$ and assume for a contradiction that $a\wedge_Q q$ and $a\vee_Q q$ do not exist.

By the Claim, $a\vee s$ does not exist, so $a\wedge s$ exists by $\BW(P,s)$.

If $a<r$, then $a<r<q$ so $a\wedge q$ would exist and equal $a\in Q$, and hence $a\wedge_Q q$ would exist, a contradiction.  Hence $r<a$.

As $r\vee s$ exists, by $\BW(P,s)$ $a\wedge(r\vee s)=a\wedge q$ exists.  Note that $r\le a\wedge q$.  If $r<a\wedge q$, then $a\wedge q\in Q$ so $a\wedge q=a\wedge_Q q$.  Thus $r=a\wedge q$.

As $s\le r\vee s=q$; $r\nleqslant s$; and $a\wedge s\le s$, we have $a\wedge s<a\wedge q=r$.  Hence $r$ has a lower cover.  By $\BWI(P,s,r)$, $r$ has at least two distinct lower covers, $r_1$ and $r_2$.  By Lemma 11, $r_1,r_2<s$ and $r_1,r_2<r<a$; hence $r_1,r_2\le a\wedge s<r$, implying that $r_1=a\wedge s=r_2$, a contradiction.

We have proven that for all $a\in Q$, $a\vee_Q q$ or $a\wedge_Q q$ exists.

Now let $a,b\in Q$ be such that $a>b$.  Assume that $b\vee_Q q$ exists but not $a\vee_Q q$.  By the Claim, $a\vee s$ does not exist.

If $a<r$, then by Lemma 11, $a<s$, so $a\vee s$ would be $s$, a contradiction.  Hence $a>r$.

{\it Case 1. $b<r$}

Then $b<r<q$, so $b\vee q=q=b\vee_Q q$.  We need to show that $a\wedge_Q q$ exists.

But it does, since $a\vee_Q q$ does not exist.

{\it Case 2. $b>r$}

We will show that $b\vee_Q q=b\vee s$: Now $b\vee_Q q\in P$ and $b\vee_Q q\ge b$ and $b\vee_Q q\ge q=r\vee s\ge s$.  If $w\in\{b,s\}^{\it u}$, then, since $r<b$, we have $w\in\{r,s\}^{\it u}$, so $w\ge r\vee s=q$.  Also, since $w\ge b>r$, we have $w\in Q$.  Thus $w\ge b\vee_Q q$.  This proves that $b\vee_Q q=b\vee s$.

By $\BW(P,s)$, $a\wedge(b\vee s)$ exists.

Now $r<a$ and $r<b\le b\vee s$, so $r\le a\wedge(b\vee s)=a\wedge(b\vee_Q q)$.  If $r<a\wedge(b\vee_Q q)$, we would be done, since $a\wedge(b\vee_Q q)$ would then be in $Q$ and hence equal to $a\wedge_Q(b\vee_Q q)$.

So assume for a contradiction that $r=a\wedge(b\vee s)$.  Then $b=a\wedge b\le a\wedge(b\vee s)=r<b$, a contradiction. $\qed$

\begin{theorem} Let $\BW(P,s)$ hold. Then $\Delta(P)$ is non-evasive.
\end{theorem}

\proof We prove this by induction on $|P|$.

Corollary 10 deals with the case $P\setminus(\uparrow s\cup\downarrow s)=\emptyset$, which includes the base case.

Now assume $P\setminus(\uparrow s\cup\downarrow s)\ne\emptyset$.  Then $\BW(P,s,r)$ holds for some $r$.  By Lemma 12 and the induction hypothesis, $\Delta(P\setminus\{r\})=\dl_{\Delta(P)}(r)$ is non-evasive.

If $r$ has a unique lower cover, then, by Lemma 5, $\Delta(P)$ is non-evasive.

So now assume that $\BWI(P,s,r)$ holds.  By Lemma 13 and the induction hypothesis, $\Delta\big((\uparrow r\cup\downarrow r)\setminus\{r\})$ is non-evasive, but this is $\lk_{\Delta(P)}(r)$.  Hence $\Delta(P)$ is non-evasive. $\qed$

\begin{corollary} Let $P$ be a finite poset with an element $s$ such that
\begin{itemize} \item[(1)] for all $x\in P$, either $s\vee x$ or $s\wedge x$ exists; and
\item[(2)] for all $x,y\in P$ such that $x<y$, if $s\wedge x$ does not exist but $s\wedge y$ does exist, then $(s\wedge y)\vee x$ exists.
\end{itemize}

Then $\Delta(P)$ is non-evasive.
\end{corollary}

\proof This is the dual of Theorem 14. $\qed$

\medskip
\medskip
\medskip

The posets in this paper arise in a still-open problem: Bj\"orner conjectures on \cite[p. 98]{BjoHA} ``that for every noncomplemented lattice $L$ of finite length the subposet $L-\{\hat0,\hat1\}$ has the fixed point property.'' (A poset $P$ has the {\it fixed point property} if, for every $f\in P^P$, there exists $p\in P$ such that $f(p)=p$.)  This was proven topologically for $L$ finite \cite[\S3]{BacBjoGI}. (Obviously we assume $|L|>2$!)

While Baclawski made a great advance in \cite{BacAB}, finding a non-topological proof for finite $L$, it is unclear how to generalize that proof to the finite {\sl length} case.  (The stumbling block is figuring out how to extend to infinite lattices the graph theory argument on page 1012 of the proof of \cite[Theorem 32]{BacAB}.)

At the 1981 Banff Conference on Ordered Sets, Bj\"orner conjectured that the {\sl finite length} analogues of the posets in Corollary 15 have the fixed point property \cite[8.8 Conjecture, p. 838]{RivHB}.

Could topological methods (perhaps replacing dismantlability with the notion of collapsibility) be extended to finite length posets?  (See, for ideas, \cite{LiMilIB} or \cite{ZhoIH}; admittedly we have not read the latter, as we still seek a translation of the latter into English.)



\begin{thebibliography}{99}

\bibitem{BacAB} Kenneth Baclawski, ``A Combinatorial Proof of a Fixed Point Property,'' {\it Journal of Combinatorial Theory (A)} {\bf 119} (2012), 994--1013.

\bibitem{BacBjoGI} Kenneth Baclawski and Anders Bj\"orner, ``Fixed Points in Partially Ordered Sets,'' {\it Advances in Mathematics} {\bf 31} (1979), 263--287.

\bibitem{BjoHA} Anders Bj\"orner, ``Homotopy Type of Posets and Lattice Complementation,'' {\it Journal of Combinatorial Theory (A)} {\bf 30} (1981), 90--100.

\bibitem{BjoIE} Anders Bj\"orner, ``Topological Methods,'' in R. L. Graham, M. Gr\"otschel, and L. Lov\'asz (eds.), {\it Handbook of Combinatorics, Volume II} (Elsevier Science B. V., Amsterdam, 1995), 1819--1872.

\bibitem{BjoWalHC} Anders Bj\"orner and James W. Walker, ``A Homotopy Complementation Formula for Partially Ordered Sets,'' {\it European Journal of Combinatorics} {\bf 4} (1983), 11--19.

\bibitem{DavPriJB} B. A. Davey and H. A. Priestley, {\it Introduction to Lattices and Order} (Cambridge University Press, Cambridge, 2002), second edition.

\bibitem{FelLueNesJG} Stefan Felsner, Marco L\"ubbecke, and Jarik Ne\v{s}et\v{r}il, ``Editorial,'' {\it European Journal of Combinatorics} {\bf 28} (2007), 2053--2056.

\bibitem{KahSakStuHD} Jeff Kahn, Michael Saks, and Dean Sturtevant, ``A Topological Approach to Evasiveness,'' {\it Combinatorica} {\bf 4} (1984), 297--306.

\bibitem{KozIH} Dmitry N. Kozlov, ``Order Complexes of Noncomplemented Lattices Are Nonevasive,'' {\it Proceedings of the American Mathematical Society} {\bf 126} (1998), 3461--3465.

\bibitem{LiMilIB} Boyu Li and E. C. Milner, ``The PT-order and the Fixed Point Property,'' {\it Order} {\bf 9} (1992), 321--331.

\bibitem{RivHB} Ivan Rival (ed.), {\it Ordered Sets: Proceedings of the NATO Advanced Study Institute Held at Banff, Canada, August 28 to September 12, 1981} (D. Reidel Publishing Company, Dordrecht, Holland, 1982).

\bibitem{SchAF} Bernd S. W. Schr\"oder, {\it Ordered Sets: An Introduction with Connections from Combinatorics to Topology} (Birkh\"auser Verlag, 2016), second edition.

\bibitem{ZhoIH} Zhou Caijun, ``The Acyclicity and Contractibility of Posets,'' {\it Acta Mathematica Sinica} {\bf 41} (1998), 361--364.



\end{thebibliography}
\end{document}